\documentclass[11pt]{article}
\usepackage[fleqn]{amsmath}
\usepackage[fleqn]{mathtools}
\usepackage[utf8]{inputenc}
\usepackage[english]{babel}
\usepackage[T1]{fontenc}
\usepackage{mathtools}
\usepackage{graphicx}
\usepackage{tabularx}
\usepackage{tikz}
\usepackage{bm}
\usepackage{bbm}

\usepackage{appendix}
\usepackage[round]{natbib}
\usetikzlibrary{patterns, backgrounds, calc, trees}
\usepackage{amssymb,amsmath,amsfonts,eurosym,geometry,ulem,graphicx,caption,color,setspace,sectsty,comment,footmisc,caption,natbib,pdflscape,subfigure,array,hyperref}

\usepackage{float}
\usepackage{mathtools}
\usepackage{mdframed}
\usepackage{lipsum}
\newtheorem{theorem}{Theorem}[section]
\newtheorem{lemma}{Lemma}[section]
\newtheorem{corollary}[theorem]{Corollary}

\newtheorem{assumption}{Assumption}[section]

\newcolumntype{L}[1]{>{\raggedright\let\newline\\arraybackslash\hspace{0pt}}m{#1}}
\newcolumntype{C}[1]{>{\centering\let\newline\\arraybackslash\hspace{0pt}}m{#1}}
\newcolumntype{R}[1]{>{\raggedleft\let\newline\\arraybackslash\hspace{0pt}}m{#1}}

\newmdtheoremenv{theo}{Theorem}
\usepackage{amsmath, amssymb, amsfonts}
\usepackage{amsmath}
\usepackage{listings}
\usepackage{eurosym}

\usepackage{setspace}

\usepackage{lmodern}
\usepackage{listings}
\usepackage{indentfirst}

\DeclarePairedDelimiter\norm{\lvert \lvert}{\rvert \rvert}
\usepackage{float}

\newcommand{\argmin}[1]{\underset{#1}{\operatorname{arg}\,\operatorname{min}}\;}

\usepackage[linesnumbered,ruled,vlined]{algorithm2e}
\graphicspath{ {images/} }\graphicspath{ {images/} }

\usepackage{geometry}
\begin{document}
\begin{titlepage}

\RestyleAlgo{boxruled}
\LinesNumbered
\title{\Large \bf Inference robust to outliers with $\boldsymbol{\ell_1}$-norm penalization\thanks{I thank my PhD supervisor Professor Eric Gautier for his availability and great help. I am also grateful to Anne Ruiz-Gazen, Jean-Pierre Florens, Thierry Magnac and Nour Meddahi for useful comments. I acknowledge financial support from the ERC POEMH 337665 grant.}}
\author{
{\large Jad B\textsc{eyhum}}
\footnote{jad.beyhum@gmail.com}\\\texttt{\small Toulouse School of Economics, Universit\'e Toulouse Capitole}}
 
\date{}

\maketitle

\begin{abstract}
\noindent This paper considers the problem of inference in a linear regression model with outliers where the number of outliers can grow with sample size but their proportion goes to $0$. We apply the square-root lasso estimator penalizing the $\ell_1$-norm of a random vector which is non-zero for outliers. We derive rates of convergence and asymptotic normality. Our estimator has the same asymptotic variance as the OLS estimator in the standard linear model. This enables to build tests and confidence sets in the usual and simple manner. The proposed procedure is also computationally advantageous as it amounts to solving  a convex optimization program. Overall, the suggested approach constitutes a practical robust alternative to the ordinary least squares estimator.  \\
\vspace{0in}
\vspace{0in} 

\bigskip
\end{abstract}
\setcounter{page}{0}
\thispagestyle{empty}
\end{titlepage}
\pagebreak \newpage

\section{Introduction}
 This paper considers a linear regression model with outliers. The statistican observes a dataset of $n$  i.i.d. realizations of an outcome scalar random variables $y_i$ and a random vector of covariates $x_i$ with support in $\mathbb{R}^K$, such that $\Sigma=\mathbb{E}[x_ix_i^\top]$ is positive definite. We place ourselves in the Huber's contamination framework, that is the distribution of $(y_i,x_i)$ is a mixture between two distributions. With probability $1/2< 1-p\le 1$, it corresponds to a linear regression model with conditionally homoscedastic errors, that is there exists $\beta\in\mathbb{R}^K$ and scalar i.i.d. random variables $\epsilon_i$ such that $y_i=x_i'\beta+\epsilon_i$, $\mathbb{E}[x_i\epsilon_i]=\mathbb{E}[\epsilon_i]=0$ and $0<\text{var}[\epsilon_i^2\vert x_i]=\sigma^2<\infty$. With probability $p$, the distribution is unspecified. An observation $(y_i,x_i)$ is called an outlier when it was generated according to this unspecified distribution $G$. The goal of the statistician is to estimate the parameter $\beta$. This model can be rewritten as
\begin{equation}\label{model}\begin{array}{cc}y_i=x_i^{\top}\beta+\alpha_i +\epsilon_i& \forall i=1,\dots,n,\end{array}\end{equation}
where $\alpha_i$ is scalar random variable which is equal to $0$ when an observation is not an outlier and which dependence with $x_i$ and $\epsilon_i$ is left unrestricted. The probability that $\alpha_i$ is different from $0$ is hence, $p= \mathbb{P}\left(\alpha_i\ne 0\right)=\mathbb{E}\left[\left|\left|\alpha \right|\right|_0/n\right]$. We derive estimation results in an asymptotic where $p$ goes to $0$ as a function of the sample size $n$.

This model can represent various situations of practical interest. 
First, the statistician could be interested in $\beta$ because it corresponds to the slope of the best linear predictor of $y_i$ given $x_i$ for the observations for which $\alpha_i=0$. These coefficients are of interest because, in the presence of outliers, the slope of the best linear predictor of $y_i$ given $x_i$ for the whole population may differ greatly from $\beta$ and hence a statistical analysis based on the whole population may lead to a poor prediction accuracy for the large part of the population that are not outliers. 

Second, if $\beta$ is given a causal interpretation, then it represents the causal effect of the regressors for the population of "standard" individuals. That is, for instance, if the aim is evaluate a program, it could be that the treatment effect is negative for most of the population but strongly positive for a small fraction of the individuals, the outliers. The policy maker may not be willing to implement a policy that has a negative effect on most of the population, giving interest to a statistical procedure that estimates the treatment effect of the majority of the population robustly.

 Finally, $\beta$ could represent the true coefficient of the best linear predictor of $y_i$ given $x_i$ in a measurement errors model. Indeed, assume that our population follows the model $\tilde{y}_i=\tilde{x}_i\beta +\tilde{\epsilon}_i$ with $\mathbb{E}[\tilde{x}_i\tilde{\epsilon}_i]=0$ but that we do not observe $(\tilde{y}_i,\tilde{x}_i)$ but $(y_i,x_i)$, this fits in our framework with $\epsilon_i=\tilde{\epsilon}_i$ and 
$$\alpha_i= y_i-\tilde{y}_i+(\tilde{x}_i-x_i)\beta.$$
Hence, $\alpha_i$ allows for both measurement errors in $x_i$ - called outliers in the $x$-direction - and in $y_i$, the outliers in the $y$-direction, for a small fraction of the population, see \citet{rousseeuw2005robust} for a precise discussion.

This paper develop results on the estimation of $\beta$ when the vector $\alpha=\left(\alpha_1,\dots,\alpha_n\right)^\top$ is sparse in the sense that $p$ goes to $0$ with $n$. We rely on a variant of the square-root lasso estimator of \citet{belloni2011square} which penalizes the $\ell_1$-norm of the vector $\alpha$. The advantages of our estimator are that the penalty parameter does not depend on the variance of the error term and is computationally tractable. If the vector $\alpha$ is sparse enough, we show that our estimator is $\sqrt{n}$-consistent and asymptotically normal. It has the same asymptotic variance as the OLS estimator in the standard linear model without outliers. 

\textbf{Related literature.} This paper is connected to at least two different research fields. First, it draws on the literature on inference in the high-dimensional linear regression model and closely related variants of this model. A series of papers from \citet{belloni2011?1,belloni2012sparse,belloni2014high,belloni2014inference,belloni2016inference,belloni2017program}  study a variety of models ranging from panel data models to quantile regression in an high-dimensional setting. \citet{gautier2011high} proposes inference procedures in an high-dimensional IV model with a large number of both regressors and instrumental variables. \citet{javanmard2014confidence, van2014asymptotically, zhang2014confidence} suggest debiasing strategies of the lasso estimator to obtain confidence intervals in a high-dimensional linear regression model. We borrow from this literature by using an $\ell_1$-penalized estimator and complete existing research by deriving inference results for the linear regression model with outliers.

Next, our work is related to the literature on robust regression. For detailed accounts of this field, see \citet{rousseeuw2005robust,hampel2011robust,maronna2018robust}. The literature identifies a trade-off between efficiency and robustness, as explicited below. Indeed, $M$-estimators (such as the Ordinary Least-Squares (OLS) estimator) are efficient when data is generated by the standard linear model without outliers and Gaussian errors but this comes at the price of a breakdown point - the maximum proportion of the data that can be contaminated without the estimator performing arbitrarily poorly - of $0$. By contrast, $S$-estimators such as the Least Median of Squares (LMS) and the Least Trimmed Squares (LTS) have a strictly positive and fixed breakdown point. They are also asymptotically normal in the model without outliers but are not efficient and have computational issues because of the non-convexity of their objective functions (see \citet{rousseeuw2005robust}). Our estimator is efficient under certain conditions, because it attains the same asymptotic variance as the OLS estimator in the standard linear model. Unlike this literature, our procedure relies on a convex program and is computationally tractable, see \citet{belloni2011square} for a detailed analysis. The proposed approach therefore provides a simple efficient alternative to the rest of the literature.

Within the robust regression literature some authors have considered the application of $\ell_1$-norm penalization to robust estimation. In particular, our model nests the Huber's contamination model for location estimation introduced in \citet{huber1964robust}. Indeed, if there is a single constant regressor, our model nests the following framework:
$$y_i=\beta+\alpha_i+\epsilon_i,$$
where $\epsilon_i\sim \mathcal{N}(0,1)$  i.i.d., $\beta\in\mathbb{R}$ is the mean of $y_i$ for non-outlying coefficients while $\mathbb{E}[y_i\vert\alpha_i\ne 0]$ is left unrestricted. \citet{chen2018robust} show that the minimax lower bound for the squared $\ell_2$-norm estimation error is of order greater than $\max(1/n,p^2)$ under gaussian errors, where $ \norm{\alpha}_0$ is the number of outliers in the sample. When $p\sqrt{\log(n)}\to 0$, we attain this lower bound up to a factor $\log(n)^2$. Several strategies have been proposed to tackle this location estimation problem. The one which is the closest to ours is soft-thresholding using a lasso estimator, that is use 
$$\widehat{\beta}\in \argmin{\beta \in \mathbb{R}^K} \sum_{i=1}^n (y_i-\beta-\alpha_i)^2 +\lambda \sum_{i=1}^n|\alpha_i|,\ \lambda>0,$$
see for instance \citet{collier2017rate}. We substitute this estimator with a square-root lasso that has the advantage to provide guidance on the choice of the penalty level that is independent from the variance of the noise (see \citet{belloni2011square}). We extend the analysis of this type of estimators to the linear regression model and add inference results to the literature. Other $\ell_1$-norm penalized estimators for robust linear regression have been studied in the literature such as in \citet{lambert2011robust,dalalyan2012socp,li2012simultaneous,alfons2013sparse}, but the authors do not provide inference results. \citet{fan2017estimation} considers robust estimation in the case where $\beta$ is a high-dimensional parameter. Its estimator penalizes the Huber loss function by a term proportional to the $\ell_1$-norm of $\beta$.

\textbf{Notations.} We use the following notations. For a matrix $M$, $M^{\top}$ is its transpose, $\norm{M}_2$ is its $\ell_2$-norm, $\norm{M}_1$ is the $\ell_1$-norm, $\norm{M}_{\infty}$ is its sup-norm, $\norm{M}_{\text{op}}$ is its operator norm and $\norm{M}_0$ is the number of non-zero coefficients in $M$, that is its $\ell_0$-norm. For a probabilistic event $\mathcal{E}$, the fact that it happens w.p.a. $1$ (with probability approaching $1$) signifies that $\mathbb{P}\left(\mathcal{E}\right)\xrightarrow[n\to \infty]{} 1$. Then, for $k=1,\dots,K$, $x_k$ is the vector $((x_{1})_k,\dots,(x_{n})_k)^\top$ and $X$ is the matrix $(x_1,\dots,x_n)^\top$. $P_X$ is the projector on the vector space spanned by the columns of the matrix $X$ and $M_X=I_n-P_X$, where $I_n$ is the identity matrix of size $n$. We introduce $y=(y_1,\dots,y_n)^\top$ and $\epsilon=(\epsilon_1,\dots,\epsilon_n)^\top$.
\section{Low-dimensional linear regression}
\label{sec:2}
\subsection{Framework}

The probabilistic framework consists of a sequence of data generating processes (henceforth, DGPs) that depend on the sample size $n$. The joint distribution of $(x_i,\epsilon_i)$ is independent from the sample size. We consider an asymptotic where $n$ goes to $\infty$ and where $p$, the contamination level, depends on $n$ while the number of regressors remains fixed.

Our estimation strategy is able to handle models where $\alpha$ is  sparse, that is $\left|\left|\alpha\right|\right|_0/n=o_P(1)$ or, in other words, $p\to 0$. Potentially, every individual's $y_i$ can be generated by a distribution that does not follow a linear model but the difference between the distribution of $y_i$ and the one yielded by a linear model can only be important for a negligible proportion of individuals. Our subsequent theorems will help to quantify these previous statements.

\subsection{Estimation procedure}

We consider an estimation procedure that estimates both the coefficients $\alpha_i$ and the effects of the regressors $\beta$ by a square-root lasso that penalizes only the coefficients $\alpha_i$, that is
$$(\widehat{\beta},\widehat{\alpha})\in \argmin{\beta\in \mathbb{R}^K,\ \alpha\in \mathbb{R}^n} \frac{1}{\sqrt{n}} \norm{y-X\beta-\alpha}_2 +\frac{\lambda}{n} \norm{\alpha}_1, $$
where $\lambda$ is a penalty level whose choice is discussed later. The advantage of the square-root lasso over the lasso estimator is that the penalty level does not depend on an estimate of the variance of $\epsilon_i$. Hence, our procedure is simple in that it does not make use of any tuning parameter unlike the least median of quares and least trimmed squares estimators. An important remark is that if $\beta$ is such that $X\beta = P_X(y-\widehat{\alpha})$, then
$$\frac{1}{\sqrt{n}} \norm{y-X\beta-\widehat{\alpha}}_2 +\frac{\lambda}{n} \norm{\widehat{\alpha}}_1\le \frac{1}{\sqrt{n}} \norm{y-Xb-\widehat{\alpha}}_2 +\frac{\lambda}{n} \norm{\widehat{\alpha}}_1,$$
for any $b\in \mathbb{R}^K$.
Therefore, if $X^\top X$ is positive definite, $\widehat{\beta}$ is the OLS estimator of the regression of $y-\widehat{\alpha}$ on $X$, that is
\begin{equation}\label{OLS}\widehat{\beta}=\left(X^\top X\right)^{-1}X^\top (y-\widehat{\alpha}).\end{equation}
Then, notice also that for all $\alpha\in \mathbb{R}^n$ and $b \in \mathbb{R}^K$, we have 
$$\frac{1}{\sqrt{n}} \norm{M_X(y-\alpha)}_2 +\frac{\lambda}{n} \norm{\alpha}_1\le \frac{1}{\sqrt{n}} \norm{y-Xb-\alpha}_2 +\frac{\lambda}{n} \norm{\alpha}_1 .$$
Hence, because $\frac{1}{\sqrt{n}} \norm{M_X(y-\alpha)}_2 +\frac{\lambda}{n} \norm{\alpha}_1$ is feasible, it holds that
\begin{equation} \label{progalpha} \widehat{\alpha}\in \argmin{\alpha\in \mathbb{R}^N} \frac{1}{\sqrt{n}}\norm{M_X(y-\alpha)}_2 +\frac{\lambda}{n} \norm{\alpha}_1.\end{equation}
Under assumptions developed below, this procedure yields consistent estimation and asymptotic normality for  $\widehat{\beta}$. Remark that model \eqref{model} can be seen as a standard linear model with the coefficient $\alpha_i$ corresponding to the slope parameter of a dummy variable which value is $1$ for the individual $i$ and $0$ otherwise. Hence, our analysis of the square-root lasso fits in the framework of \citet{belloni2011square}. However, our approach is met with additional technical difficulties because we penalize only a subset of the variables and there is no hope to estimate $\alpha$ consistently as each of its entries is indirectly observed only once.  As a result, we develop new assumptions and theorems that are better suited for the purposes of this paper.

\subsection{Assumptions and results}
The main assumption concerns the choice of the penalty level:
\begin{assumption}
\label{as1.1}
We have $\lim \limits_{n\to \infty}\mathbb{P}\left(\lambda\ge 2\sqrt{n}\frac{\norm{M_X\epsilon}_{\infty}}{\norm{M_X\epsilon}_2}\right)=1.$
\end{assumption}
The tuning of $\lambda$ prescribed by this penalty level depends on the distributional assumptions made on $\epsilon$, in particular on the tails. The next lemma provides guidance on how to choose the regularization parameter according to assumptions on $\epsilon$: 
\begin{lemma}
\label{choice}
It holds that $2\sqrt{n}\frac{\norm{M_X\epsilon}_{\infty}}{\norm{M_X\epsilon}_2}=2\frac{\norm{\epsilon}_{\infty}}{\sigma}+o_P(\norm{\epsilon}_{\infty})+O_P(1)$. Additionally, if $\psi$ is such that $\lim \limits_{n\to \infty}\mathbb{P}\left(\psi\ge 2\frac{\norm{\epsilon}_{\infty}}{\sigma}\right)=1$ and $\varphi\to \infty$, then for any $c>1$, $\lambda=c\psi+\varphi$ satisfies Assumption \ref{as1.1}.
\end{lemma}
The proof is given in Appendix A. This lemma suppresses the role of the matrix $X$ in the choice of the penalty and simplifies the decision procedure. It leads to the subsequent corollary:
\begin{corollary}The following hold:
\label{sub}\begin{itemize}
\item[(i)] If $\epsilon_i$ are gaussian random variables, then $\lambda=2c\sqrt{2\log(n)}$ satisfies Assumption \ref{as1.1} for any $c>1$;
\item[(ii)] If $\epsilon_i$ are sub-gaussian random variables, then there exists a constant $c>0$ such that $\lambda=c\sqrt{\log(n)}$ satisfies Assumption \ref{as1.1};
\item[(iii)] If $\epsilon_i$ are sub-exponential random variables, then  then there exists a constant $c>0$ such that $\lambda=c\log(n)$ satisfies Assumption \ref{as1.1}.
\end{itemize}
\end{corollary}
The proof is given in Appendix A. The statistician can use Corollary \ref{sub} to decide on the penalization parameter given how heavy she expects the tails of the error term to be in her data. In practice, it is advised to choose the smallest penalty verifying Assumption \ref{as1.1}. This can be done by Monte-Carlo simulations. Notice that our approach allows for heavy-tailed distributions such as sub-exponential random variables.

To derive the convergence rate of our estimator, we first bound the estimation error on $\alpha$ and obtain the following result:
\begin{lemma}\label{alpha} Under Assumption \ref{as1.1} and if $p\max\left(\lambda,\sqrt{\left|X\right|_{\infty}}\right)=o_P(1)$ (and , it holds that 
 $$\frac{1}{n}\norm{\widehat{\alpha}-\alpha}_1 =O_{P}\left( p\lambda\right).$$
 \end{lemma}
The proof is given in Appendix B. The rate of convergence of $\norm{\widehat{\alpha}-\alpha}_1/n$ therefore is lower than $p\sqrt{\log(n)}$ if the errors are gaussian or sub-gaussian and we choose the penalty level as in Lemma \ref{sub}.
Note that, as standard in works related to the lasso estimator (see \citet{buhlmann2011statistics}), in our proof we make use of a compatibility condition that states that a compatibility constant is bounded from below with probability approaching one. The condition that $p\norm{X}_{\infty}=o_P(1)$ is enough to show that this property holds as shown in Lemma \ref{lA.3} in Appendix B. It is possible to find other sufficient conditions but it is outside the scope of this paper. Remark that if $\{x_i\}_i$ are i.i.d. sub-Gaussian random variables then $\norm{X}_{\infty}=O_P\left(\sqrt{\log(n)}\right)$ allowing for the sparsity level $p=o_P(1/\sqrt{\log(n)})$. 

Here, we show how to derive the rate of convergence of $\widehat{\beta}$ from Lemma \ref{alpha}. Assume that $p\max\left(\lambda,\left|X\right|_{\infty}\right)=o_P(1)$. Substituting $y$ by $X\beta+\alpha+\epsilon$ in \eqref{OLS}, we obtain
\begin{equation} \label{decompOLS} \widehat{\beta}-\beta=\left(X^{\top}X\right)^{-1}X^{\top}\epsilon +(X^{\top}X)^{-1}X^{\top}(\alpha-\widehat{\alpha}).\end{equation}
Now, notice that $\left(X^{\top}X\right)^{-1}X^{\top}(\alpha-\widehat{\alpha})=(X^{\top}X/n)^{-1}X^{\top}(\alpha-\widehat{\alpha})/n$. By the law of large numbers, $\left(X^{\top}X/n\right)^{-1}=O_{P}(1)$, which implies that 
\begin{align}\notag \left|\left|\left(X^{\top}X\right)^{-1}X^{\top}(\alpha-\widehat{\alpha})\right|\right|_2&\le\left|\left\vert \left(\frac{1}{n}X^{\top}X\right)^{-1}\right\vert \right\vert_{\text{op}}\frac{1}{n}\left|\left|X^{\top}(\alpha-\widehat{\alpha})\right|\right|_{2}\\
&=O_{P}\left(\frac{1}{n}\norm{X}_{\infty}\norm{\alpha-\widehat{\alpha}}_1\right)\quad \text{(By H\"older's inequality)} \label{influential}.
\end{align}
  By Lemma \ref{alpha}, this implies that $$\norm{(X^{\top}X)^{-1}X^{\top}(\alpha-\widehat{\alpha})}_2=O_{P}\left(p\lambda\norm{X}_{\infty}\right).$$
 Finally, by the central limit theorem and Slutsky's lemma, we have that $\sqrt{n}(X^{\top}X)^{-1}X^{\top}\epsilon  \xrightarrow[]{d}\mathcal{N}(0,\sigma \Sigma^{-1})$. This leads to Theorem \ref{th1.1}.
 \begin{theorem}\label{th1.1} Under Assumption \ref{as1.1} and if $p\max\left(\lambda,\sqrt{\left|X\right|_{\infty}}\right)=o_P(1)$, it holds that 
   $$\frac{\widehat{\beta}-\beta}{\max\left( \frac{1}{\sqrt{n}},p\lambda\norm{X}_{\infty}\right)}=O_P\left(1\right).$$
   \end{theorem}
  
This result allows to derive the rates of convergence under different assumptions on the tails of the distributions of the regressors and the error term. For instance, if $\{x_i\}_i$ and $\{\epsilon_i\}_i$ are i.i.d. sub-Gaussian random variables, then $\widehat{\beta}$ is consistent as long as $p\log(n)\to 0$ for the choice of $\lambda$ proposed in Lemma \ref{sub}. In this case, this implies that our estimator reaches (up to a logarithmic factor) the minimax lower bound for the Huber's contamination location model under gaussian errors, which is $\max(1/n,p^2)$ in $\ell_2$-norm  according  to \citet{chen2018robust}. We attain the rate $ \max(1/n,p^2\log(n))$. Remark also that equation \eqref{influential} explains the role of $\norm{X}_{\infty}$ in the convergence rate of $\widehat{\beta}$. For an individual $i$, if $x_i$ is large then an error in the estimation of $\alpha_i$ can contribute to an error in the estimation of $\beta$ via the term $(X^{\top}X)^{-1}X^{\top}(\alpha-\widehat{\alpha})$ in \eqref{decompOLS}. $\norm{X}_{\infty}$ measures the maximum influence that an observation can have.

To show that our estimator is asymptotically normal, it suffices to assume that the bias term $(X^{\top}X)^{-1}X^{\top}(\alpha-\widehat{\alpha})$ in \eqref{decompOLS} vanishes asymptotically: 
\begin{theorem} \label{th1.2} Under Assumption \ref{as1.1}, assuming that $p\lambda\norm{X}_{\infty}\sqrt{n}=o_P(1)$ (and $p\max\left(\lambda,\sqrt{\left|X\right|_{\infty}}\right)=o_P(1)$), we have 
   $$\sqrt{n}(\widehat{\beta}-\beta) \xrightarrow[]{d}\mathcal{N}(0,\sigma^2 \Sigma^{-1}).$$
Moreover, $\widehat{\sigma}^2=\frac{1}{n}\sum_{i=1}^n (y_i-x_i^{\top}\hat{\beta}-\widehat{\alpha})^2$ and $\widehat{\Sigma}=\frac{1}{ n}\sum_{i=1}^n x_i x_i^{\top}$ are consistent estimators of, respectively, $\sigma^2$ and $\Sigma$.
   \end{theorem}
The proof that $\widehat{\sigma}^2\xrightarrow{\mathbb{P}} \sigma^2$ is given in Appendix C.When the entries of $X$ and $\epsilon$ are sub-Gaussian, for the choice of the penalty prescribed in Lemma \ref{sub}, the contamination level needs to satisfy $p\log(n)\sqrt{n}\to 0$ to be able to use \ref{th1.2} to prove asymptotic normality. Notice that the asymptotic variance of our estimator corresponds to the one of the OLS estimator in the standard linear model under homoscedasticity. Hence, confidence sets and tests can be built in the same manner as in the theory of the OLS estimator. 

An important last remark concerns the meaning of confidence intervals developed using Theorem \ref{th1.2}. Note that they are obtained under an asymptotic with triangular array data under which the number of outliers is allowed to go to infinity while the proportion of outliers goes to $0$. The interpretation of a 95\% confidence interval $I$ built with Theorem \ref{th1.2} is as follows: if the number of outliers in our data is low enough and the sample size is large enough, then there is a probability of approximatively $0.95$ that $\beta$ belongs to $I$.

\section{Computation and simulations}
\subsection{Iterative algorithm}
We propose the following algorithm to compute our estimator. Because $u=\min_{\sigma>0}\left\{\frac{\sigma}{2}+\frac{1}{2\sigma}u^2\right\}$, as long as $\left|\left\vert y-X\widehat{\beta}-\widehat{\alpha}\right|\right\vert_2^2>0 $, we have that
\begin{equation}\label{optiglob}
(\widehat{\beta},\widehat{\alpha}, \widehat{\sigma})\in \arg\min_{\beta\in \mathbb{R}^K,\alpha\in \mathbb{R}^n,\sigma\in \mathbb{R}^+}\frac{\sigma}{2}+ \frac{1}{2\sigma}
\left|\left\vert y-X\beta-\alpha \right|\right\vert_2^2+\frac{\lambda}{2\sqrt{n}}\left|\left|\alpha\right|\right|_1.
\end{equation}
This is a convex objective and we propose to iteratively minimize over $\beta$, $\alpha$, and $\sigma$. Let us start from $\left(\beta^{(0)},\alpha^{(0)},\sigma^{(0)}\right)$ and compute the following sequence for $t\in\mathbb{N}^*$ until convergence:
\begin{enumerate}
\item $\beta^{(t+1)}\in \arg\min_{\beta\in \mathbb{R}^K} \left|\left\vert y-X\beta-\alpha^{(t)}\right\vert\right|_2^2;$
\item $\alpha^{(t+1)}\in \arg \min_{\alpha\in \mathbb{R}^n}  \left\vert\left| y-X\beta^{(t+1)}-\alpha \right|\right\vert_2^2+\frac{\lambda\sigma^{(t)}}{\sqrt{n}}\left|\left|\alpha\right|\right|_1;$
\item $\sigma^{(t+1)}=\left| \left\vert y-X\beta^{(t+1)}-\alpha^{(t+1)}\right|\right\vert_2$.
\end{enumerate}
The following lemma explains how to perform step 2:
\begin{lemma}\label{alphaiter}
For $i=1,\dots,n$, if $\left|y_i-(X\beta^{(t+1)})_i \right|\le \frac{\lambda\sigma^{(t)}}{\sqrt{n}}$ then $\alpha^{(t+1)}_i=0$. If $\left|y_i-(X\beta^{(t+1)})_i \right|> \frac{\lambda\sigma^{(t)}}{\sqrt{n}}$ then $\alpha^{(t+1)}_i=y_i-\left(X\beta^{(t+1)}\right)_i -\text{sign}\left(y_i-\left(X\beta^{(t+1)}\right)_i\right)\frac{\lambda\sigma^{(t)}}{\sqrt{n}}$.
\end{lemma}
The proof is given in Appendix D.
\subsection{Simulations}
We apply this computation approach in a small simulation exercise. The data generating process is as follows: there are two regressors $x_{1i}$ and $x_{2i}$, with $x_{1i}=1$ for all $i$ and $x_{2i}$ are i.i.d. $\mathcal{N}(0,1)$ random variables. $\epsilon_i$ are i.i.d. $\mathcal{N}(0,1)$ random variables. Then, we set 
$$ \alpha_i=\left\{\begin{array}{cc} 0&\text{if $x_{2i}<q_{1-p}$}\\
5x_{2i}&\text{if $x_{2i}\ge q_{1-p}$},
\end{array}\right.$$
where $q_{1-p}$ is such that $\mathbb{P}(x_{2i}\ge q_{1-p})=p$. In table~\ref{fig:Covoutliers}, we present the bias, the variance, the mean squared error (MSE) and the coverage of $95\%$ confidence intervals for our estimator $\widehat{\beta}$ computed using the algorithm of Section 3.1, where we use $100$ iterations and with $\lambda=2.01\sqrt{2\log(n)}$. This choice corresponds to the one outlined in Corollary \ref{sub}. The bias, the variance and the coverage of $95\%$ confidence intervals for the naive OLS estimator:
$$\widetilde{\beta}^{OLS}\in \argmin{\beta\in \mathbb{R}^K} \norm{y-X\beta-\alpha}_2^2$$
are also reported.
For the OLS estimator, the confidence intervals correspond to the ones of the standard linear model. The presented results are averages among $8000$ replications. We observe that our estimator brings a substantial improvement in estimation precision with respect to the OLS estimator.

\begin{table}[!ht]
\centering
             \begin{tabular}{|l|l|l|c|c|c|c|}
  \hline
   value &p&n& $\widehat{\beta}_1$ & $\widetilde{\beta}_1^{OLS}$ &$\widehat{\beta}_2$  & $\widetilde{\beta}_2^{OLS}$ \\
  \hline
  \hline
    \text{bias}& 0.025 & 100& 0.127 & 0.301& 0.278&0.671\\
    \hline 
    \text{variance} & 0.025& 100& 0.060& 0.130& 0.097&0.221\\
    \hline 
        \text{MSE}& 0.025 & 100& 0.076 & 0.221& 0.174&0.671\\
    \hline 
    \text{coverage} &0.025 &100 & 0.82& 0.47& 0.75&0.20\\
    \hline
      \hline
    \text{bias}& 0.01 & 1000& 0.045 & 0.120& 0.133&0.361\\
    \hline 
    \text{variance} & 0.01& 1000& 0.002& 0.004& 0.007&0.018\\
    \hline 
        \text{MSE}& 0.001 & 1000& 0.004 & 0.018& 0.025&0.148\\
    \hline 
    \text{coverage} &0.001 &1000 & 0.74& 0.16& 0.24&0.00\\
    \hline
      \hline
    \text{bias}& 0.001 & 10000& 0.005 & 0.015& 0.017&$0.057$\\
    \hline 
    \text{variance} & 0.001& 10000& 1.08 $\times 10^{-4}$& 1.28 $\times 10^{-4}$& 2.21 $\times 10^{-4}$&5.23 $\times 10^{-4}$\\
    \hline 
        \text{MSE}& 0.001 & 10000& 1.33 $\times 10^{-4}$ & 3.53 $\times 10^{-4}$& 5.10 $\times 10^{-4}$&3.772 $\times 10^{-3}$\\
    \hline 
    \text{coverage} &0.001 &10000 & 0.93& 0.66& 0.68&0.03\\
    \hline
      \end{tabular}  
        \captionsetup{labelsep=none}
    \caption{. bias, variance, mean squared error (MSE) and coverage of $95\%$ confidence intervals for $\lambda=2.01\sqrt{2\log(n)}$.}
        \label{fig:Covoutliers}
    \end{table}

\bibliographystyle{plainnat}
\bibliography{refOutliers}
\appendix
\appendixpage
\section{Choice of the penalization parameter}
\subsection{Proof of Lemma \ref{choice}}
We start by proving the next two technical lemmas:
\begin{lemma}\label{tech1}It holds that $\norm{P_X\epsilon}_\infty= O_P(1)$.\end{lemma}
{\bf Proof.}  
Because of the assumptions on the joint distribution of $(x_i,\epsilon_i)$, we have that $\sqrt{n}(X^{\top}X)^{-1}X^{\top}\epsilon  \xrightarrow[]{d}\mathcal{N}(0,\sigma \Sigma^{-1})$, therefore $\sqrt{n}\norm{(X^{\top}X)^{-1}X^{\top}\epsilon}_{2}=O_{P}(1)$.
Because $X(X^{\top}X)^{-1}X^{\top}\epsilon=\frac{X}{\sqrt{n}}\sqrt{n}(X^{\top}X)^{-1}X^{\top}\epsilon$, we obtain that $\norm{P_X\epsilon}_2\le\frac{\norm{X}_2}{\sqrt{n}}\sqrt{n}\norm{(X^{\top}X)^{-1}X^{\top}\epsilon}_{2} =O_{P}\left(\frac{\norm{X}_{2}}{\sqrt{n}}\right)=O_P(1)$, by the law of large numbers. 
 \hfill $\Box$\\
\begin{lemma}\label{tech2}It holds that $\frac{\sqrt{n}}{\norm{M_X\epsilon}_2}-\frac{1}{\sigma}= o_{P}(1)$. \end{lemma}
{\bf Proof.}  
First, remark that, by the theorem of Pythagore, \begin{align}\notag \norm{M_X\epsilon}_2^2&=\left<\epsilon-X(X^{\top}X)^{-1}X^{\top}\epsilon,\epsilon-X(X^{\top}X)^{-1}X^{\top}\epsilon\right>\\ \notag & =\norm{\epsilon}_2^2-\epsilon^{\top}X(X^{\top}X)^{-1}X^{\top}\epsilon.\end{align}
Now, this leads to $\frac{1}{n}\norm{M_X\epsilon}_2^2=\frac{1}{n}\norm{\epsilon}_2^2-\frac{1}{n}\epsilon^{\top}X(X^{\top}X)^{-1}X^{\top}\epsilon$. Because $\{x_i\}_i$ and $\{\epsilon_i\}_i$ are i.i.d. and $\mathbb{E}[x_i\epsilon_i]=0$, we have that $\sqrt{n}(X^{\top}X)^{-1}X^{\top}\epsilon  \xrightarrow[]{d}\mathcal{N}(0,\sigma \Sigma^{-1})$ and $\frac{1}{\sqrt{n}}X^{\top}\epsilon  \xrightarrow[]{d}\mathcal{N}(0,\sigma \Sigma)$. This implies that $\epsilon^{\top}X(X^{\top}X)^{-1}X^{\top}\epsilon=O_{P}(1/n)$. We also have that $\frac{1}{n}\norm{\epsilon}_2^2\xrightarrow[]{\mathbb{P}}\sigma^2$, which leads to $\frac{1}{n}\norm{M_X\epsilon}_2^2\xrightarrow[]{\mathbb{P}}\sigma^2$. We conclude by the continuous mapping theorem.
 \hfill $\Box$\\

Now, we proceed with the proof of Lemma \ref{choice}. Notice that 
\begin{align*}
2\sqrt{n}\frac{\norm{M_X(\epsilon)}_{\infty}}{\norm{M_X\epsilon}_2}&\le \frac{2\sqrt{n}}{\norm{M_X\epsilon}_2}(\norm{\epsilon}_{\infty}+\norm{P_X\epsilon}_{\infty})\\
&\le \frac{1}{\sigma}\norm{\epsilon}_{\infty} + \left|\frac{\sqrt{n}}{\norm{M_X\epsilon}_2}-\frac{1}{\sigma}\right|\norm{\epsilon}_{\infty}+\frac{1}{\sigma}\norm{P_X\epsilon}_{\infty} + \left|\frac{\sqrt{n}}{\norm{M_X\epsilon}_2}-\frac{1}{\sigma}\right|\norm{P_X\epsilon}_{\infty}
\end{align*}
Using lemmas \ref{tech1} and \ref{tech2}, we obtain
\begin{equation} \label{oproba}2\sqrt{n}\frac{\norm{M_X\epsilon}_{\infty}}{\norm{M_X\epsilon}_2}=2\frac{\norm{\epsilon}_{\infty}}{\sigma}+o_P(\norm{\epsilon}_{\infty})+O_P(1)\end{equation}
The rest of the lemma is a direct consequence of \eqref{oproba} and the pigeonhole principle.
\subsection{Proof of Corollary \ref{sub}}
\noindent \textbf{Proof of (i)} By Lemma \ref{choice} it is sufficient to show that for $c>1$,
$$\lim \limits_{n\to \infty}\mathbb{P}\left(2c\sqrt{2\log(n)}\ge 2\frac{\norm{\epsilon}_{\infty}}{\sigma}\right)=1.$$
Let us remember the gaussian bound (see Lemma B.1 in \citet{giraud2014introduction}): for $t\ge 0$, we have 
$$\mathbb{P}\left(\frac{\left|\epsilon_i\right|}{\sigma}\ge t\right)\le 2 e^{-\frac{t^2}{2}}.$$
Then, we have 
\begin{align*}
\mathbb{P}\left(2c\sqrt{2\log(n)}\ge 2\frac{\norm{\epsilon}_{\infty}}{\sigma}\right)&\le \sum_{i=1}^n  \mathbb{P}\left(c\sqrt{2\log(n)}\ge \frac{\left|\epsilon_i\right|}{\sigma}\right)\\
& \le ne^{-c\log(n)}\to 0.
\end{align*}
\noindent \textbf{Proof of (ii)} By Lemma \ref{choice} it is sufficient to show that there exists $c>0$ such that
$$\lim \limits_{n\to \infty}\mathbb{P}\left(c\sqrt{\log(n)}\ge 2\frac{\norm{\epsilon}_{\infty}}{\sigma}\right)=1.$$
Let us remember the sub-gaussian bound (see Proposition 2.5.2 in \citet{vershynin2018high}): for $t\ge 0$, there exists $b>0$ such that
$$\mathbb{P}\left(\frac{\left|\epsilon_i\right|}{\sigma}\ge t\right)\le 2 e^{-\frac{t^2}{2b}}.$$
Then, we have 
\begin{align*}
\mathbb{P}\left(4\sqrt{b}\sqrt{\log(n)}\ge 2\frac{\norm{\epsilon}_{\infty}}{\sigma}\right)&\le \sum_{i=1}^n  \mathbb{P}\left(2\sqrt{b}\sqrt{\log(n)}\ge \frac{\left|\epsilon_i\right|}{\sigma}\right)\\
& \le 2ne^{-2\log(n)}\to 0.
\end{align*}
\noindent \textbf{Proof of (iii)} By Lemma \ref{choice} it is sufficient to show that there exists $c>0$ such that
$$\lim \limits_{n\to \infty}\mathbb{P}\left(c\log(n)\ge 2\frac{\norm{\epsilon}_{\infty}}{\sigma}\right)=1.$$
Let us remember the sub-exponential bound (see Proposition 2.7.1 in \citet{vershynin2018high}): for $t\ge 0$, there exists $b>0$ such that
$$\mathbb{P}\left(\frac{\left|\epsilon_i\right|}{\sigma}\ge t\right)\le 2 e^{-\frac{t}{2b}}.$$
Then, we have, for $n$ large enough,
\begin{align*}
\mathbb{P}\left(8b\sqrt{\log(n)}\ge 2\frac{\norm{\epsilon}_{\infty}}{\sigma}\right)&\le \sum_{i=1}^n  \mathbb{P}\left(4b\sqrt{\log(n)}\ge \frac{\left|\epsilon_i\right|}{\sigma}\right)\\
& \le 2ne^{-2\log(n)}\to 0.
\end{align*}
\section{Proof of lemma \ref{alpha}}
\subsection{Compatibility constant}
For $\delta\in \mathbb{R}^n$, we denote by $\delta_{J}\in \mathbb{R}^n$ the vector for which $(\delta_{J})_i=\delta_i$ if $\alpha_i\ne 0$ and $(\delta_{J})_i=0$ otherwise. Let us also define $\delta_{J^c}=\delta-\delta_{J}$.
We introduce the following cone:
$$C=\left\{\delta\in \mathbb{R}^n \ s.t.\ \left|\left|\delta_{J^c}\right|\right|_1 \le 3\left|\left|\delta_{J}\right|\right|_1 \right\}.$$
We work with the following compatibility constant (see \citet{buhlmann2011statistics} for a discussion of the role of compatibility conditions in the lasso literature) corresponding to
$$\kappa =\min_{\delta\in C, \delta \ne 0}\frac{\sqrt{2\norm{\alpha}_0}\norm{M_X\delta}_{2}}{\left|\left|\delta_{J}\right|\right|_1}.$$
We use the following lemma:
  \begin{lemma} \label{lA.3}
If $p^2\norm{X}_{\infty}=o_P(1)$, there exists $\kappa_*>0$ such that $\kappa>\kappa_*$ w.p.a. $1$.
\end{lemma}
 {\bf Proof.}  
 Take $\delta \in C$, to show this result, notice that
$$M_X\delta=\delta -X(X^{\top}X)^{-1}X^{\top}\delta.$$
Therefore, we have 
\begin{align}
\notag\norm{M_X\delta}_2&\ge \norm{\delta}_2-\norm{X(X^{\top}X)^{-1}X^{\top}\delta}_2\\
\notag &= \norm{\delta}_2-\left|\left|\sum_{k=1}^KX_k\left((X^{\top}X)^{-1}X^{\top}\delta\right)_k\right\vert\right\vert_2\\
\notag &\ge \norm{\delta}_2-\sum_{k=1}^K\left|\left|X_k\left((X^{\top}X)^{-1}X^{\top}\delta\right)_k\right\vert\right\vert_2\\
\notag &\ge \norm{\delta}_2-\sum_{k=1}^K\left|\left|X_k\right|\right|_2\left|\left|(X^{\top}X)^{-1}X^{\top}\delta\right\vert\right\vert_{\infty}\\
\notag &\ge \norm{\delta}_2-\sum_{k=1}^K\left|\left|X_k\right|\right|_2\left|\left|(X^{\top}X)^{-1}X^{\top}\delta\right\vert\right\vert_{2}\\
\notag &\ge \norm{\delta}_2-\sum_{k=1}^K\norm{X_k}_{2}\left\vert \left\vert \left(\frac{1}{n}X^{\top}X\right)^{-1}\right\vert \right\vert_{\text{op}} \frac{1}{n}\norm{X^{\top}\delta}_2\\
\notag &\ge \norm{\delta}_2-\sum_{k=1}^K\norm{X_k}_{2}\left\vert \left\vert \left(\frac{1}{n}X^{\top}X\right)^{-1}\right\vert \right\vert_{\text{op}}\frac{\sqrt{K}}{n} \norm{X}_{\infty}\norm{\delta}_1 \quad \text{(By H\"older's inequality)}\\
\notag &\ge \norm{\delta}_2-\sum_{k=1}^K\norm{X_k}_{2}\left\vert \left\vert \left(\frac{1}{n}X^{\top}X\right)^{-1}\right\vert \right\vert_{\text{op}}\frac{\sqrt{K}}{n} \norm{X}_{\infty}4\norm{\delta_J}_1 \quad \text{(Because $\delta \in C$)}\\
\notag &\ge \norm{\delta}_2-\sum_{k=1}^K\norm{X_k}_{2}\left\vert \left\vert \left(\frac{1}{n}X^{\top}X\right)^{-1}\right\vert \right\vert_{\text{op}}\frac{\sqrt{K}}{n} \norm{X}_{\infty}4\sqrt{\norm{\alpha}_0}\norm{\delta_J}_2 \quad \text{(Because $\norm{\delta_J}_0\le \norm{\alpha}_0$)}\\
\label{derniereligne} &\ge \norm{\delta}_2-\sum_{k=1}^K\frac{\norm{X_k}_{2}}{\sqrt{n}}\left\vert \left\vert \left(\frac{1}{n}X^{\top}X\right)^{-1}\right\vert \right\vert_{\text{op}} 4\sqrt{K}\sqrt{\frac{\norm{\alpha}_0}{n}}\norm{X}_{\infty}\norm{\delta}_2, \end{align}
where $X_k=(x_{1k},...,x_{nK})^{\top}$.
Next, we have that 
\begin{align*}
\kappa &\ge \min_{\delta\in C, \delta \ne 0}\frac{\sqrt{2\norm{\alpha}_0}\norm{M_X\delta}_{2}}{\left|\left|\delta_{J}\right|\right|_1}\\
&\ge \min_{\delta\in C, \delta \ne 0}\frac{\sqrt{2\norm{\alpha}_0}\norm{M_X\delta}_{2}}{\sqrt{||\alpha||_0}\left|\left|\delta_{J}\right|\right|_2}\\
&\ge \sqrt{2} \min_{\delta\in C, \delta \ne 0}\frac{\norm{M_X\delta}_{2}}{\left|\left|\delta\right|\right|_2}\\
&\ge \sqrt{2} \left(1-\sum_{k=1}^K\frac{\norm{X_k}_{2}}{\sqrt{n}}\left\vert \left\vert \left(\frac{1}{n}X^{\top}X\right)^{-1}\right\vert \right\vert_{\text{op}} 4\sqrt{K}\sqrt{\frac{\norm{\alpha}_0}{n}}\norm{X}_{\infty}\right)
\end{align*}
Now, because we have $\frac{1}{n}\sum_{i=1}^n x_{i}x_i^{\top}\xrightarrow[]{\mathbb{P}}\Sigma$ by the law of large numbers, we obtain that $\left\vert \left\vert \left(X^{\top}X/n\right)^{-1}\right\vert \right\vert_{\text{op}}=O_{P}(1)$ and that $\sum_{k=1}^K\norm{X_k}_2/\sqrt{n}=\sum_{k=1}^K\sqrt{\left(X^{\top}X/n\right)_{kk}}=O_{P}(1)$, both implying that$\frac{1}{\sqrt{n}}\sum_{k=1}^K\norm{X_k}_2\left\vert \left\vert \left(X^{\top}X/n\right)^{-1}\right\vert \right\vert_{\text{op}}=O_{P}(1)$. We conclude the proof using that $p^2\norm{X}_{\infty}=o_P(1)$.  \hfill $\Box$\\

\subsection{End of the proof of Lemma \ref{alpha}}
Throughout this proof, we work on the event 
$$\left\{\lambda \ge \frac{2\sqrt{n}\norm{M_X\epsilon}_2{\infty}}{\norm{M_X\epsilon}_2}\right\}\cap\left\{\kappa>\kappa_*\right\}\cap \left\{\left(\frac{2\sqrt{2p}\lambda}{\kappa}\right)^2<1\right\} ,$$ 
which has probability approaching $1$ according to Assumption \ref{as1.1}, Lemma \ref{lA.3}, and the condition that $p\lambda \to 0$. Let us define $\Delta=\widehat{\alpha}-\alpha$. Now, remark that
\begin{align}
 \notag \norm{\widehat{\alpha}}_1&=\norm{\alpha+\Delta}_1\\
\notag&=\left|\left|\alpha+\Delta_{J}+\Delta_{J^c}\right|\right|_1\\
 &\ge \left|\left|\alpha+\Delta_{J^c}\right|\right|_1-\left|\left|\Delta_{J}\right|\right|_1
\label{finally}.
\end{align}
 Next, we use the fact that $\left|\left|\alpha+\Delta_{J^c}\right|\right|_1=\norm{\alpha}_1+\left|\left|\Delta_{J^c}\right|\right|_1$. Combining this and \eqref{finally}, we get
\begin{equation}\label{Trois}\norm{\widehat{\alpha}}_1\ge \norm{\alpha}_1+\left|\left|\Delta_{J^c}\right|\right|_1-\left|\left|\Delta_{J}\right|\right|_1.
\end{equation}
By definition of $\widehat{\alpha}$ and concentrating our objective function in $\beta$, we have
 \begin{equation}\label{minimizer} \frac{1}{\sqrt{n}}\norm{M_X(y-\widehat{\alpha})}_2+\frac{\lambda}{n} \norm{\widehat{\alpha}}_1\le \frac{1}{\sqrt{n}}\norm{M_X(y-\alpha)}_2+\frac{\lambda}{n} \norm{\alpha}_1.\end{equation}
 By convexity, if $M_X\epsilon \ne 0$, it holds that
     \begin{align}
     \notag \frac{1}{\sqrt{n}}\norm{M_X(y-\widehat{\alpha})}_2-\frac{1}{\sqrt{n}}\norm{M_X(y-\alpha)}_2&\ge -\frac{1}{\sqrt{n}\norm{M_X\epsilon}_2}\left<M_X(\epsilon),\Delta\right>\\
    \label{ineq}  &\ge-\frac{\lambda}{2n} \norm{\Delta}_1,
     \end{align}
 where \eqref{ineq} comes from $\lambda \ge 2\sqrt{n}\norm{M_X\epsilon}_2/\norm{M_X\epsilon}_{\infty}$. This last inequality is also straightforwardly true when $M_X\epsilon=0$. 
 This and \eqref{minimizer} imply
  \begin{equation}\label{biz}\norm{\hat{\alpha}}_1\le \frac{1}{2} \norm{\Delta}_1+\norm{\alpha}_1.\end{equation}
  Using \eqref{Trois}, we get
  $$\norm{\alpha}_1+\left|\left|\Delta_{J^c}\right|\right|_1-\left|\left|\Delta_{J}\right|\right|_1\le \frac{1}{2} \norm{\Delta}_1+\norm{\alpha}_1.$$
Then, as $\norm{\Delta}_1=\left|\left|\Delta_{J^c}\right|\right|_1+\left|\left|\Delta_{J}\right|\right|_1$, we obtain
   \begin{equation}
   \label{when}
 \left|\left|\Delta_{J^c}\right|\right|_1\le 3 \left|\left|\Delta_{J}\right|\right|_1,\end{equation}
which implies that $\Delta\in C$. 
Using $y=X\beta+\alpha+\epsilon$, we get
   $$\frac{1}{n}\norm{M_X(y-\widehat{\alpha})}_2^2-\frac{1}{n}\norm{M_X(y-\alpha)}_2^2=\frac{1}{n}\norm{M_X(\widehat{\alpha}-\alpha)}_2^2-\frac{2}{n}\left<M_X\epsilon,\widehat{\alpha}-\alpha\right>.$$
By H\"older's inequality, this results in 
$$\frac{1}{n}\norm{M_X(y-\widehat{\alpha})}_2^2-\frac{1}{n}\norm{M_X(y-\alpha)}_2^2\le \frac{1}{n}\norm{M_X(\widehat{\alpha}-\alpha)}_2^2-\frac{2}{n}\norm{M_X\epsilon}_{\infty}\norm{\Delta}_1.$$
     Because $\lambda \ge 2\sqrt{n}\frac{\norm{M_X\epsilon}_{\infty}}{{\norm{M_X\epsilon}_2}}$ , we obtain
     $$\frac{1}{n}\norm{M_X(\widehat{\alpha}-\alpha)}_2^2\le \frac{1}{n}\norm{M_X(y-\widehat{\alpha})}_2^2-\frac{1}{n}\norm{M_X(y-\alpha)}_2^2+\frac{\lambda\norm{M_X\epsilon}_2}{n^{\frac{3}{2}}}\norm{\Delta}_1.$$
This implies that
     \begin{align}
     \notag &\frac{1}{n}\norm{M_X(\widehat{\alpha}-\alpha)}_2^2&\\
     \notag & \le \frac{1}{n}\norm{M_X(y-\widehat{\alpha})}_2^2-\frac{1}{n}\norm{M_X(y-\alpha)}_2^2+\frac{\lambda\norm{M_X\epsilon}_2}{n^{\frac{3}{2}}}\norm{\Delta}_1\\
     \notag &  = \frac{1}{n}\norm{M_X(y-\widehat{\alpha})}_2^2-\frac{1}{n}\norm{M_X(y-\alpha)}_2^2+\frac{\lambda\norm{M_X\epsilon}_2}{n^{\frac{3}{2}}}\left(\left|\left|\Delta_{J}\right|\right|_1+\left|\left|\Delta_{J^c}\right|\right|_1\right)\\
      &   \label{One}\le \frac{1}{n}\norm{M_X(y-\widehat{\alpha})}_2^2-\frac{1}{n}\norm{M_X(y-\alpha)}_2^2+\frac{4\lambda\norm{M_X\epsilon}_2}{n^{\frac{3}{2}}}\left|\left|\Delta_{J}\right|\right|_1\quad \text{(Because $\Delta \in C$)}.
     \end{align}
By equations \eqref{Trois} and \eqref{minimizer}, we have $ \frac{1}{\sqrt{n}}\norm{M_X(y-\widehat{\alpha})}_2-\frac{1}{\sqrt{n}}\norm{M_X(y-\alpha)}_2\le \frac{\lambda}{n}\left(\left|\left|\Delta_{J}\right|\right|_1-\left|\left|\Delta_{J^c}\right|\right|_1\right)$. Using the fact that $\Delta \in C$ and \eqref{ineq}, this yields $$\left \vert  \frac{1}{\sqrt{n}}\norm{M_X(y-\widehat{\alpha})}_2-\frac{1}{\sqrt{n}}\norm{M_X(y-\alpha)}_2\right \vert \le\frac{2\lambda}{n}\left|\left|\Delta_{J}\right|\right|_1.$$
Next, notice that  \begin{align*}&\frac{1}{n}\norm{M_X(y-\widehat{\alpha})}_2^2-\frac{1}{n}\norm{M_X(y-\alpha)}_2^2\\
&=\left(\frac{1}{\sqrt{n}}\norm{M_X(y-\widehat{\alpha})}_2-\frac{1}{\sqrt{n}}\norm{M_X(y-\alpha)}_2\right)\left(\frac{1}{\sqrt{n}}\norm{M_X(y-\widehat{\alpha})}_2+\frac{1}{\sqrt{n}}\norm{M_X(y-\alpha)}_2\right).\end{align*}
This implies
   \begin{align}
   \notag &\left\vert \frac{1}{n}\norm{M_X(y-\widehat{\alpha})}_2^2-\frac{1}{n}\norm{M_X(y-\alpha)}_2^2 \right \vert\\
  \notag  &\le \frac{2\lambda}{n}\left|\left|\Delta_{J}\right|\right|_1\left( \frac{2}{\sqrt{n}}\norm{M_X(y-\alpha)}_2+\frac{2\lambda}{n}\left|\left|\Delta_{J}\right|\right|_1\right)\\
   \label{Two}  &\le \left(\frac{2\lambda}{n}\right)^2\left|\left|\Delta_{J}\right|\right|_1^2+ \frac{4}{\sqrt{n}}\norm{M_X(y-\alpha)}_2\frac{\lambda}{n}\left|\left|\Delta_{J}\right|\right|_1.
   \end{align}
Combining \eqref{One} and \eqref{Two} and remarking that $\norm{M_X\epsilon}_2=\norm{M_X(y-\alpha)}_2$, we obtain
 $$ \frac{1}{n}\norm{M_X(\widehat{\alpha}-\alpha)}_2^2\le  \left(\frac{2\lambda}{n}\right)^2\left|\left|\Delta_{J}\right|\right|_1^2+\frac{4\norm{M_X\epsilon}_2}{\sqrt{n}}\frac{\lambda}{n}\left|\left|\Delta_{J}\right|\right|_1+\frac{4\lambda\norm{M_X\epsilon}_2}{n^{\frac{3}{2}}}\left|\left|\Delta_{J}\right|\right|_1.$$
 Now, as $\Delta \in C$, this implies that
 $$ \frac{1}{n}\norm{M_X\Delta}_2^2\le  \left(\frac{2\lambda}{n}\right)^2\left(\frac{\sqrt{2\norm{\alpha}_0}\norm{M_X\Delta}_{2}}{\kappa}\right)^2+\frac{8\lambda\norm{M_X\epsilon}_2}{n^{\frac{3}{2}}}\frac{\sqrt{2\norm{\alpha}_0}\norm{M_X\Delta}_{2}}{\kappa}.$$
From now on assume that $\norm{M_X\Delta}_2\ne 0$, we get
 $$\frac{1}{n}\norm{M_X\Delta}_2\le \left(1- \left(\frac{2\sqrt{2\frac{\norm{\alpha}_0}n}\lambda}{\kappa}\right)^2\right)^{-1}\frac{8\norm{M_X\epsilon}_2\sqrt{2\frac{\norm{\alpha}_0}n}\lambda}{n\kappa},$$
  which implies again that 
$$\frac{1}{n}\left|\left|\Delta_{J}\right|\right|_1\le \left(1- \left(\frac{2\sqrt{2\frac{\norm{\alpha}_0}n}\lambda}{\kappa}\right)^2\right)^{-1}\frac{16\norm{M_X\epsilon}_2\frac{\norm{\alpha}_0}n\lambda}{\sqrt{n}\kappa^2}. $$
Finally, as $\Delta \in C$, we have
 \begin{align}
 \notag \frac{1}{n}\norm{\Delta}_1&=\frac1n\left(\left|\left|\Delta_{J}\right|\right|_1+\left|\left|\Delta_{J^c}\right|\right|_1\right)\\
 \notag &\le \frac1n4\left|\left|\Delta_{J}\right|\right|_1\\
 \label{last} &\le \left(1- \left(\frac{2\sqrt{2\frac{\norm{\alpha}_0}n}\lambda}{\kappa_*}\right)^2\right)^{-1}\frac{64\norm{M_X\epsilon}_2\frac{\norm{\alpha}_0}n\lambda}{\sqrt{n}\kappa_*}.
 \end{align}
The last inequality also holds if $M_X\Delta=0$ because, as $\kappa>\kappa_*$, this implies that $\Delta_J=0$ and hence $\Delta=0$ using the fact that $\Delta $ belongs to $C$. To conclude the proof, use \eqref{last}, the fact that $\norm{M_X\epsilon}_2/\sqrt{n}\le \norm{\epsilon}_2/\sqrt{n}=o_P(1)$ by the law of large numbers and the condition $p\max\left(\lambda,\norm{X}_{\infty}\right)=o_P(1)$.
\section{Proof that \textbf{$\mathbf{\widehat{\sigma}^2\xrightarrow{\mathbb{P}} \sigma^2}$} in Theorem \ref{th1.2}}
We have 
\begin{align*}
\widehat{\sigma}^2&=\frac{\left|\left|y-X\widehat{\beta}-\widehat{\alpha}\right|\right|_2^2}{n}\\
&= \frac{\left|\left|X\left(\beta-\widehat{\beta}\right)+\left(\alpha-\widehat{\alpha}\right)+\epsilon\right|\right|_2^2}{n}\\
&= \frac{\left|\left|X\left(\beta-\widehat{\beta}\right)+\left(\alpha-\widehat{\alpha}\right)\right|\right|_2^2}{n}+2\frac{\left<X\left(\beta-\widehat{\beta}\right)+\left(\alpha-\widehat{\alpha}\right),\epsilon\right>}{n}+2\frac{\left|\left|\epsilon\right|\right|_2^2}{n}\\
&= \frac{\left|\left|X\left(\beta-\widehat{\beta}\right)\right|\right|_2^2}{n}+\frac{\left|\left|\alpha-\widehat{\alpha}\right|\right|_2^2}{n}+2\frac{\left<X\left(\beta-\widehat{\beta}\right),\alpha-\widehat{\alpha}\right>}{n}\\
&\quad+2\frac{\left<X\left(\beta-\widehat{\beta}\right),\epsilon\right>}{n}+2\frac{\left<\alpha-\widehat{\alpha},\epsilon\right>}{n}+2\frac{\left|\left|\epsilon\right|\right|_2^2}{n}.
\end{align*}
Next, remark that there exists $c,\eta>0$ such that $\mathbb{P}\left(\left|X\right|_{\infty}\ge c\right)\ge \delta$ for $n$ large enough. Indeed otherwise, this would imply that the law of large numbers cannot hold for $x_i$. This yields that $\sqrt{n}p\lambda=o(1)$. Then, because of Lemma \ref{alpha}, Theorem \ref{th1.2} and $p\lambda\left|X\right|_{\infty}=o_P(1)$ it holds that 
\begin{align*}
\left|\left|\widehat{\alpha}-\alpha\right|\right|_1&=o_P\left(\sqrt{n}\right)
;\\
\left|\left|\widehat{\beta}-\beta\right|\right|_2&=o_P\left(\frac{1}{\sqrt{n}}\right).
\end{align*}
Next, we have 
\begin{align*}
\frac{\left|\left|X\left(\beta-\widehat{\beta}\right)\right|\right|_2^2}{n}&\le \frac{\left|\left|X\right|\right|_2^2\left|\left|\widehat{\beta}-\beta\right|\right|_2^2}{n}\\
&=o_P(1),
\end{align*}
by the law of large numbers.
Then, by H\"older's inequality, we obtain that 
\begin{align*}
\frac{\left|\left|\alpha-\widehat{\alpha}\right|\right|_2^2}{n}&\le \frac{\left|\left|\alpha-\widehat{\alpha}\right|\right|_1\left|\left|\alpha-\widehat{\alpha}\right|\right|_\infty}{n}\\
&\le \frac{\left|\left|\alpha-\widehat{\alpha}\right|\right|_1^2}{n}\\
&=o_P(1).
\end{align*}
By the inequality of Cauchy-Schwartz, this also leads to $\frac{\left<X\left(\beta-\widehat{\beta}\right),\alpha-\widehat{\alpha}\right>}{n}=o_P(1)$. Then, the law of large numbers implies that $\frac{\left<X\left(\beta-\widehat{\beta}\right),\epsilon\right>}{n}=o_P(1)$ and $\frac{\left<\alpha-\widehat{\alpha},\epsilon\right>}{n}=o_P(1)$, which concludes the proof.
\section{Proof of Lemma \ref{alphaiter}}
By Lemma D.5 in \citet{giraud2014introduction}, there exists $\widehat{z}$, a random vector in $\mathbb{R}^n$, such that the first-order conditions of step 2 are
\begin{equation}\label{FOCMARCH} -\left(y-X\beta^{(t)}-\alpha^{(t+1)}\right)+ \frac{\lambda\sigma^{(t)}}{\sqrt{n}}\widehat{z}=0,\end{equation}
where, for $i=1,\dots,n$, $\widehat{z}_i\in [-1,1]$ if $\alpha^{(t+1)}_i=0$ and $\widehat{z}_i=\text{sign}\left(\alpha^{(t+1)}_i\right)$ if $\alpha^{(t+1)}_i\ne0$. This yields that, if $\alpha^{(t+1)}_i\ne0$, 
$$ \alpha^{(t+1)}_i=y_i-\left(X\beta^{(t+1)}\right)_i -\text{sign}\left(\alpha^{(t+1)}_i\right)\frac{\lambda\sigma^{(t)}}{\sqrt{n}}.$$
Hence, if $\alpha^{(t+1)}_i>0$, we obtain
$$ \alpha^{(t+1)}_i=y_i-\left(X\beta^{(t+1)}\right)_i -\frac{\lambda\sigma^{(t)}}{\sqrt{n}}$$
and, therefore, $y_i-\left(X\beta^{(t+1)}\right)_i >\frac{\lambda\sigma^{(t)}}{\sqrt{n}}\ge 0$. Similarly, if $\alpha^{(t+1)}_i<0$, we have
$$ \alpha^{(t+1)}_i=y_i-\left(X\beta^{(t+1)}\right)_i +\frac{\lambda\sigma^{(t)}}{\sqrt{n}}$$
and, therefore, $y_i-\left(X\beta^{(t+1)}\right)_i <-\frac{\lambda\sigma^{(t)}}{\sqrt{n}}\le 0$. This shows that, if $\alpha^{(t+1)}_i\ne 0$, we have
$$\alpha^{(t+1)}_i=y_i-\left(X\beta^{(t+1)}\right)_i -\text{sign}\left(y_i-\left(X\beta^{(t+1)}\right)_i\right)\frac{\lambda\sigma^{(t)}}{\sqrt{n}}$$ 
and $$\left|y_i-(X\beta^{(t+1)})_i \right|> \frac{\lambda\sigma^{(t)}}{\sqrt{n}}.$$ Next, if $\alpha^{(t+1)}_i=0$, \eqref{FOCMARCH} implies that
$$\left|y_i-(X\beta^{(t+1)})_i \right|\le \frac{\lambda\sigma^{(t)}}{\sqrt{n}}.$$ This concludes the proof.
 \end{document}